%% file: arxivv2.tex
\documentclass[12pt]{article}

\usepackage[paper=letterpaper,margin=.7in,twoside=false,includehead]{geometry}

\usepackage{amsfonts,amsmath,amsthm, amssymb} 
\usepackage{hyperref}
\usepackage[capitalize]{cleveref}
\usepackage{graphicx}
\usepackage{fancyhdr}
\usepackage{xparse}
\usepackage{mathtools}
\usepackage{thm-restate}
\usepackage[makeroom]{cancel}
\usepackage[shortlabels]{enumitem}
\usepackage[normalem]{ulem}

\SetEnumitemKey{thmparts}{topsep=2pt plus 1pt minus1pt, itemsep=1pt plus0.5pt minus0.5pt}

\usepackage{xcolor}
\usepackage{verbatim}

\input basics

\newtheorem{thm}{Theorem}
\newtheorem{lem}[thm]{Lemma}
\newtheorem{cor}[thm]{Corollary}
\newtheorem{conj}[thm]{Conjecture}

\theoremstyle{definition}

\newtheorem{example}[thm]{Example}

\newtheorem{claim}{Claim}
\numberwithin{claim}{thm}
\newenvironment{proofc}{\begin{proof}[Proof of Claim]}{\end{proof}}
\newtheorem{case}{Case}
\numberwithin{case}{claim}

\numberwithin{subcase}{case}
\newtheorem{prop}[thm]{Proposition}

\newcommand{\G}[2]{G_n^{#1}(#2)}
\newcommand{\GOne}[1]{\G{#1}{i_d}}
\newcommand{\GTwo}[1]{\G{#1}{i_0}}

\newcommand{\calG}[2]{\mathcal{G}_n^{#1}(#2)}

\newcommand{\h}[2]{H_n^{#1}(#2)}

\begin{document}

\title{Maximizing the number of stars in graphs with forbidden properties}
\author{
Zhanar Berikkyzy\footnote{Department of Mathematics, Fairfield University, Fairfield, CT, USA  Email: {\tt zberikkyzy@fairfield.edu}. Supported in part by NSF grant DMS-2418903.}
\qquad
Kirsten Hogenson\footnote{Department of Mathematics and Statistics, Skidmore College, Saratoga Springs, NY, USA  Email: {\tt khogenso@skidmore.edu}.}
\qquad
Rachel Kirsch\footnote{Department of Mathematical Sciences, George Mason University, Fairfax, VA, USA  Email: {\tt rkirsch4@gmu.edu}. Supported in part by Simons Foundation Grant MP-TSM-00002688.}
\qquad
Jessica McDonald\footnote{Auburn University, Department of Mathematics and Statistics, Auburn, AL, USA
  Email: {\tt mcdonald@auburn.edu}.   
	Supported in part by Simons Foundation Grant \#845698.  }
}

\date{}

\maketitle

\begin{abstract} 
Erd\H{o}s proved an upper bound on the number of edges in an $n$-vertex non-Hamiltonian graph with given minimum degree and showed sharpness via two members of a particular graph family. F\"{u}redi, Kostochka and Luo showed that these two graphs play the same role when ``number of edges'' is replaced by ``number of t-stars,'' and that two members of a more general graph family maximize the number of edges among non-$k$-edge-Hamiltonian graphs. In this paper we generalize their former result from Hamiltonicity to related properties (traceability, Hamiltonian-connectedness, $k$-edge Hamiltonicity, $k$-Hamiltonicity) and their latter result from edges to $t$-stars. We identify a family of extremal graphs for each property that is forbidden. 
This problem without the minimum degree condition was also open; here we conjecture a complete description of the extremal family for each property, and prove the characterization in some cases. Finally, using a different family of extremal graphs, we find the maximum number of $t$-stars in non-$k$-connected graphs. 

\end{abstract}

\section{Introduction}

In this paper all graphs are simple.

Let $n,i,\ell\in\mathbb{Z}$ with $-1\le\ell\leq n-3$ and $1\leq i\leq \tfrac{n-1-\ell}{2}$. Define $\G{\ell}{i}$ to be the graph $K_{i+\ell}+(I_i\cup K_{n-2i-\ell})$, where $\cup$ indicates disjoint union and $+$ indicates a complete bipartite graph between the two sets of vertices. (In the case when $i+\ell=0$, we define $\G{\ell}{i}$ to be $I_1 \cup K_{n-1}$.) Two members of this family will be of particular interest for us, namely when $i$ takes on its top value of $i_0:=\floor{\tfrac{n-1-\ell}{2}}$ and, for any nonnegative integer $d\leq \tfrac{n-1+\ell}{2}$, when $i$ takes on the value of $i_d:=\max\{1, d-\ell\}$ (note $i_d \le i_0$).  Erd\H{o}s used the $\ell=0$ case of this graph family in \cite{ErdosRemarks}, noting that $\G{0}{i}$ is a \emph{non-Hamiltonian} graph with minimum degree $i$, and proving the following theorem. Note that by a graph being \emph{Hamiltonian}, we mean that it contains a \emph{Hamilton cycle}, that is, a spanning cycle.

\begin{thm}[Erd\H{o}s \cite{ErdosRemarks}] \label{Erd} Let $G$ be an $n$-vertex graph with minimum degree $\delta(G)\geq d$, where $1\leq d\leq \lfloor \tfrac{n-1}{2}\rfloor$. If $G$ is not Hamiltonian, then $e(G)\leq \max\{e(\G{0}{i_0}), e(\G{0}{i_d})\}.$
\end{thm}

F\"{u}redi, Kostochka and Luo \cite{FKL18} found that not only do $\G{0}
{i_0}, \G{0}{i_d}$ maximize the number of edges among nonhamiltonian graphs with $n$ vertices and minimum degree at least $d$, but they also maximize the number of different $t$-stars in such a graph. Note that given a graph $G$ and a 
vertex $v\in V(G)$ with degree $d$, $G$ contains $\binom{d}{t}$ different \emph{$t$-stars centered at $v$}, that is, $\binom{d}{t}$ different copies of $K_{1, t}$ where $v$ is the 
vertex of degree $t$, for any $t\in\mathbb{Z}^+$. Let $s_t(G)$ be the number of $t$-stars in a given graph $G$, so 
$s_t(G)=\sum_{v\in V(G)}\binom{d(v)}{t}$ for $t \ge 2$ and 
$s_t(G) = e(G) = \frac{1}{2}\sum_{v \in V(G)} d(v)$ for $t=1$. 
F\"{u}redi, Kostochka, Luo \cite{FKL18} proved the following.

\begin{thm}[F\"{u}redi, Kostochka, Luo \cite{FKL18}] \label{FKL} 
Let $G$ be an $n$-vertex graph with minimum degree\\ $\delta(G)\geq d$, where $1\leq d\leq \lfloor \tfrac{n-1}{2}\rfloor$, and let $t\in\{1, \ldots, n-1\}$. If $G$ is not Hamiltonian, then 
$s_t(G) \le \max\{s_t(\G{0}{i_0}), s_t(\G{0}{i_d})\}.$
\end{thm}

In this paper we generalize Theorem \ref{FKL} from Hamiltonicity to the following related properties. 
Given a graph $G$, a \emph{Hamilton path} is a path containing every vertex of $G$; $G$ is \emph{traceable} if it contains a \emph{Hamilton path}, and \emph{Hamiltonian-connected} if it contains a \emph{Hamilton path} between every pair of distinct vertices. One way to describe how ``strongly'' Hamiltonian an $n$-vertex graph $G$ is, is to say that $G$ is \emph{$k$-edge Hamiltonian}, which means that every linear forest of size at most $k$ is contained in a Hamilton cycle of $G$, for some $k\in\{0, \ldots, n-3\}$. On the other hand, we can give a measure of ``robustness'' of Hamiltonicity for an $n$-vertex graph $G$ by saying that $G$ is \emph{$k$-Hamiltonian}, which means that the removal of any set of at most $k$ vertices results in a Hamiltonian graph, for some $k\in\{0, \ldots, n-3\}$. The properties $0$-edge Hamiltonicity and $0$-Hamiltonicity are equivalent to Hamiltonicity.

The generalized graphs $\G{\ell}{i}$ defined above have previously appeared in the literature, for example in \cite[Theorems 3.16 and 3.18]{BBHKNSWY15} and in the following result.

\begin{thm}[F\"{u}redi, Kostochka, Luo \cite{FKL19}]\label{thm:FKL2}
    Let $G$ be an $n$-vertex graph with minimum degree\\ $\delta(G)\geq d$, where $k+1\leq d\leq \lfloor \tfrac{n+k-1}{2}\rfloor$. If $G$ is not $k$-edge Hamiltonian, then 
$e(G) \le \max\{e(\G{k}{i_0}),$ $e(\G{k}{i_d})\}.$
\end{thm}

F\"{u}redi, Kostochka, and Luo also proved that the same two graphs from Theorem \ref{thm:FKL2} achieve the maximum number of $t$-cliques. We instead generalize from edges to $t$-stars.

Define $\mathcal{G}^{\ell}_n(i)$ to be the set of all spanning subgraphs of $\G{\ell}{i}$ where the only edges allowed to be missing are those from the $K_{n-2i-\ell}$.  Just as the graph $\G{0}{i}$ is a non-Hamiltonian graph with minimum degree $i$, it turns out that all graphs in the families $\calG{\ell}{i}$ have minimum degree $i+\ell$ and do not have the analogous forbidden properties. 

\begin{restatable}{prop}{properties}\label{prop:properties} Let $n,i,\ell\in\mathbb{Z}$, with $-1\le\ell\leq n-3$ and $1\leq i\leq \tfrac{n-1-\ell}{2}$. For every graph $G$ in $\calG{\ell}{i}$, $G$: is not Hamiltonian when $\ell=0$; is not traceable when $\ell=-1$; is not Hamiltonian-connected when $\ell=1$; is not $k$-edge Hamiltonian when $\ell=k$ for some  $k\in\{1, \ldots, n-3\}$; and is not $k$-Hamiltonian when $\ell=k$ for some $k\in\{1, \ldots, n-3\}$.
\end{restatable}

Our first main result is the following theorem.

\begin{thm}\label{thm:main} Let $G$ be an $n$-vertex graph with minimum degree $\delta(G)\geq d$, and let $t\in\{1,2, \ldots, n-1\}$. Suppose that at least one of the following is true:
\begin{enumerate}
\item[(H1)] $G$ is not Hamiltonian;
\item[(H2)] $G$ is not traceable; 
\item[(H3)] $G$ is not Hamiltonian-connected;
\item[(H4)] $G$ is not $k$-edge Hamiltonian for some integer $1 \le k \le n-3$; or
\item[(H5)] $G$ is not $k$-Hamiltonian for some integer $1\leq k\leq n-3$.
\end{enumerate}
Then let $\ell$ be $0, -1, 1, k, k$, for the cases (H1)-(H5), respectively, and suppose that $0\leq d\leq \lfloor\tfrac{n+\ell-1}{2}\rfloor$.
Then $s_t(G) \le \max\set{s_t(\G{\ell}{i_d}),s_t(\G{\ell}{i_0})}$, and this bound is tight.
\end{thm}

Note that case (H1) of Theorem \ref{thm:main} is precisely Theorem \ref{FKL}, and the other cases are variations. Also, case (H4) when $t=1$ is Theorem \ref{thm:FKL2}. We will refer to cases (H1)--(H5) regularly throughout this paper.

In addition to cases (H1)-(H5), we get an analogue to Theorem \ref{thm:main} for $k$-connectedness. The extremal examples look a little different however. For $n, k, i \in \mathbb{Z}$, where $1 \le k \le n-2$ and $1 \le i \le \frac{n-k+1}{2}$, we define $\h{k}{i}$ to be the graph $K_{k-1} + (K_i \cup K_{n-k-i+1})$. (In the case when $k=1$, we define $\h{1}{i}$ to be $K_i\cup K_{n-i}$.) Notice that $\h{k}{i}$ is not $k$-connected because the vertices of the $K_{k-1}$ form a cut set of size less than $k$.

\begin{restatable}{thm}{connectedthm}\label{thm:kconn}
    Let $G$ be an $n$-vertex graph with $n \ge 3$ and minimum degree $\delta(G) \ge d$ for some $0 \le d \le (n+k-3)/2$, and let $t \in \set{1, \ldots, n-1}$. If $G$ is not $k$-connected, then \[s_t(G) \le 
        \begin{cases}
            s_t(\h{k}{i_d}) & \text{for } t=1\\
            \max\set{s_t(\h{k}{i_d}), s_t(\h{k}{i_0})} & \text{for }2 \le t \le n-1,
        \end{cases}
    \]
    and this bound is tight.
\end{restatable}

Returning now to our focus on cases (H1)-(H5) in Theorem \ref{thm:main},  there are two questions that arise naturally in each case, namely: (1) \emph{``Which of $s_t(\G{\ell}{i_d}),s_t(\G{\ell}{i_0})$ is the larger quantity?''} and (2) \emph{``If we knew this, could we describe the extremal family completely?''} It turns out that the second question is the easier of the two.  To this end, recall that $\G{\ell}{i}$ is the graph $K_{i+\ell}+(I_i\cup K_{n-2i-\ell})$, and $\mathcal{G}^{\ell}_n(i)$ is the set of all spanning subgraphs of $\G{\ell}{i}$ where the only edges allowed to be missing are those from the $K_{n-2i-\ell}$.  Observe that $\G{\ell}{i_0}$ has $K_{n-2i_0-\ell}=K_1$ if $n \not\equiv\ell \pmod{2}$ and $K_{n-2i_0-\ell}=K_2$ if $n \equiv\ell \pmod{2}$, so the family $\mathcal{G}^{\ell}_n(i_0)$ contains only one or two graphs. 
However, the family $\mathcal{G}^{\ell}_n(i_d)$ may have many members.

This leads to our second main theorem.

\begin{restatable}{thm}{extremalgraphs}\label{thm:extFam}
    \begin{enumerate}
    \item If $s_t(\G{\ell}{i_d}) <s_t(\G{\ell}{i_0})$ in Theorem \ref{thm:main}, then for $t \le n-i_0-1$, $\G{\ell}{i_0}$ is the unique extremal graph achieving this upper bound, and for $t > n-i_0-1$ the set of all graphs achieving this upper bound is precisely $\mathcal{G}^{\ell}_n(i_0)$. 
    \item If $s_t(\G{\ell}{i_d}) > s_t(\G{\ell}{i_0})$ in Theorem \ref{thm:main}, then for $t \le n-i_d-1$, $\G{\ell}{i_d}$ is the unique extremal graph achieving this upper bound, and for $t > n-i_d-1$, the set of all graphs achieving this upper bound is precisely $\mathcal{G}^{\ell}_n(i_d)$. 
    \item If $s_t(\G{\ell}{i_d}) = s_t(\G{\ell}{i_0})$ in Theorem \ref{thm:main}, then for $t \le n-i_0-1$, $\set{\G{\ell}{i_d},\G{\ell}{i_0}}$ is the set of extremal graphs; for $n-i_0-1 < t \le n-i_d-1$, $\mathcal{G}^{\ell}_n(i_0)\cup \set{\G{\ell}{i_d}}$ is the set of extremal graphs; and for $t > n-i_d-1$, $\mathcal{G}^{\ell}_n(i_d) \cup  \mathcal{G}^{\ell}_n(i_0)$ is the set of extremal graphs.
    \end{enumerate}
\end{restatable}

Note that in the case $\ell=0$ F\"{u}redi, Kostochka, Luo \cite[Claim 12]{FKL18} claimed that either $\G{0}{i_d}$ or $\G{0}{i_0}$ is the unique extremal graph, but the other members of $\mathcal{G}^{0}_n(i_d)$ or $\mathcal{G}^{0}_n(i_0)$ are also extremal for large values of $t$. See \cref{ex:otherextremal} (at the end of \cref{sec:mainproof}) for a concrete example.

In order to consider the difficult question of (1) above \cref{thm:extFam}---\emph{``Which of $s_t(\G{\ell}{i_d}),s_t(\G{\ell}{i_0})$ is the larger quantity?''}--- it simplifies matters to remove the minimum degree condition from Theorem \ref{thm:main}, and instead consider the following corollary.

\begin{cor}\label{cor:main} Let $G$ be an $n$-vertex graph and let $t\in\{1,2, \ldots, n-1\}$. Suppose that at least one of (H1)--(H5) is true, and set $\ell$ equal to $0, -1, 1, k, k$, respectively. Then 
$s_t(G) \le \max\set{s_t(\G{\ell}{1}),$ $s_t(\G{\ell}{i_0})}$, and this bound is tight.
\end{cor}

We prove the following, our third main theorem.

\begin{thm}\label{l0conj} 
In Corollary \ref{cor:main}, when $\ell=0$, for $n\geq 4$,
\begin{align*}
        s_t(\G{0}{1}) \le s_t(\G{0}{i_0}) & \quad\text{for } t \geq\tfrac{n+1}{2},\\
        s_t(\G{0}{i_0}) \le s_t(\G{0}{1}) & \quad\text{for } t < \tfrac{n+1}{2},
\end{align*}
and the inequalities are strict for $n \ge 6$. For all $\ell$, if $t\geq\tfrac{n+\ell+1}{2}$, then 
$s_t(\G{\ell}{1}) \le s_t(\G{\ell}{i_0})$, and the inequality is strict for $0\leq \ell \le n-5$.
\end{thm}

Note that a discussion of part of the $\ell=0$ case of Theorem \ref{l0conj} appeared in \cite{FKL18}. We conjecture that the bound of $t\geq\tfrac{n+\ell+1}{2}$ in Theorem \ref{l0conj} is tight in the sense that, for sufficiently large $n$, it is exactly the threshold for which one of the two quantities is larger. 

\begin{conj}\label{bigConj} In Corollary \ref{cor:main}, if $t < \tfrac{n+\ell+1}{2}$, then $s_t(\G{\ell}{i_0}) < s_t(\G{\ell}{1})$ for sufficiently large $n$.  
\end{conj}

Our paper now proceeds as follows. Section \ref{sec:prelims} contains some preliminaries for our work, including a proof of Proposition \ref{prop:properties}.  The proof of Theorem \ref{thm:main} appears in Section \ref{sec:mainproof} of this paper, and the proof of Theorem \ref{thm:kconn} is presented in Section \ref{sec:kconn}. We prove Theorem \ref{l0conj} in Section \ref{sec:bigConj}.

\section{Preliminaries}\label{sec:prelims}

We address five properties simultaneously using the fact that they are all \emph{$s$-stable}. 
    A property $P$ is \emph{$s$-stable} if, for all graphs $G$ and nonadjacent vertices $u$ and $v$ in $G$, whenever $G+uv$ has $P$ and $d_G(u) + d_G(v) \ge s$, the graph $G$ itself has $P$.

For example, the fact that Hamiltonicity is $n$-stable was proved by Ore \cite{ore}. It is also important that the five properties $P$ addressed in this paper all hold for sufficiently large complete graphs, so there exists an integer $n(P)$ such that $K_n$ has property $P$ for every $n \ge n(P)$. The table below shows, for each property $P$, the value of $s$ for which $P$ is $s$-stable and the value of $n(P)$. The number $\ell$ is simply $s-n$, so that each property $P$ is $(n+\ell)$-stable. 
The data in this table is presented in \cite{BC76} (using the fact that $0$-Hamiltonian-connectedness is Hamiltonian-connectedness \cite{Berge}). 

\begin{center}
\begin{tabular}{lllr}
    Property & $s$ & $n(P)$ & $\ell$\\\hline
    Traceability & $n-1$ & $2$ & $-1$\\
    Hamiltonicity & $n$ & $3$ & $0$\\
    Hamiltonian-connectedness & $n+1$ & $2$ & $1$\\
    $k$-edge Hamiltonicity & $n+k$ & $3$ & $k$\\
    $k$-Hamiltonicity & $n+k$ & $k+3$ & $k$\\
\end{tabular}
\end{center}

Bondy and Chv\'atal \cite{BC76} proved that if $P$ is $s$-stable and $n(P)$ exists then the following Chv\'atal-like degree condition holds for $P$. 
(See also \cite{BBHKNSWY15} for discussion of related properties and conditions.)

\begin{thm}[Bondy and Chv\'atal \cite{BC76}]\label{thm:Pdegree}
    Let $P$ be an $(n+\ell)$-stable property for which $n(P)$ exists. Let $G$ be an $n$-vertex graph for $n \ge n(P)$ and $d_1 \le \cdots \le d_n$ its degrees. If $G$ does not have $P$, then there is an integer $1 \le i \le \frac{n-1-\ell}{2}$ for which $d_i \le i+\ell$ and $d_{n-i-\ell} \le n-i-1$. In other words, $G$ has at least $i$ vertices of degree at most $i+\ell$ and at most $i+\ell$ vertices of degree at least $n-i$.
\end{thm}

\cref{thm:Pdegree} implies $\sigma_2$,  minimum degree, and edge conditions for these properties as well, as shown in \cite{Dawkins24, DK23}. The minimum degree and edge conditions will be relevant for this paper.

\begin{cor}[Theorem 2.3.4 in Dawkins \cite{Dawkins24}, Theorem 2.2 in Dawkins and Kirsch \cite{DK23}]\label{cor:Pdegree}
    Let $P$ be an $(n+\ell)$-stable property for which $n(P)$ exists. Let $G$ be an $n$-vertex graph for $n \ge n(P)$. If $G$ does not have $P$, then $\delta(G) \le (n+\ell-1)/2$.
\end{cor}
\begin{proof}
    Let $G$ be an $n$-vertex graph not having property $P$, with degrees $d_1 \le \cdots \le d_n$. By \cref{thm:Pdegree}, there is an integer $1 \le i \le \frac{n-1-\ell}{2}$ for which $d_i \le i+\ell$. The minimum degree of $G$ then is \[d_1 \le d_i \le i+\ell \le \frac{n-1-\ell}{2}+\ell = \frac{n+\ell-1}{2}.\qedhere\]
\end{proof}

\cref{cor:Pdegree} explains the choice to restrict $d$ to be at most $(n+\ell-1)/2$. We are looking for the maximum value of $s_t(G)$ over the set of $n$-vertex graphs that do not have $P$ and that have minimum degree at least $d$. This set of graphs would be empty if $d$ were greater than $(n+\ell-1)/2$.

The graphs $\G{\ell}{i}$ defined in the introduction are the extremal graphs for \cref{thm:Pdegree}, showing that it is best possible, as they do not have the property $P$ (\cref{prop:properties}), their degree lists are entry-wise the maximum allowed by the theorem, and they are the unique graphs for their degree lists (\cref{prop:uniquedegree}). The special case that $t=1$ and $d\le\ell+1$ (so $i_d=1$) is known:
\begin{thm}[Theorem 2.3.5 in Dawkins \cite{Dawkins24}, Theorem 2.3 in Dawkins and Kirsch \cite{DK23}]\label{thm:edgemax}
    Let $P$ be an $(n+\ell)$-stable property for which $n(P)$ exists. Let $G$ be an $n$-vertex graph with $n \ge n(P)$. If $G$ does not have $P$, then $e(G) \le e(\G{\ell}{1})$.
\end{thm}

The first piece of Theorem \ref{thm:main} we handle is that the graphs in question show the bound is tight: we show in \cref{prop:properties} that the graphs do not have the forbidden properties, and in \cref{prop:uniquedegree} that they satisfy the minimum degree condition.

\properties*

\begin{proof}  
    \begin{enumerate}
        \item Let $G \in \mathcal{G}^{0}_n(i)$. Then $G$ by definition is a spanning subgraph of $K_{i} + (I_i \cup K_{n-2i})$, where only the edges in the $K_{n-2i}$ are allowed to be missing. Deleting the $i$ dominating vertices yields a graph with at least $i+1$ components, so $G$ is not Hamiltonian.
        \item Let $G \in \mathcal{G}^{-1}_n(i)$. Then $G$ by definition is a spanning subgraph of $K_{i-1} + (I_i \cup K_{n-2i+1})$, where only the edges in the $K_{n-2i+1}$ are allowed to be missing. Deleting the $i-1$ dominating vertices yields a graph with at least $i+1$ components, so $G$ is not traceable.
        \item Let $G \in \mathcal{G}^{1}_n(i)$. Then $G$ by definition is a spanning subgraph of $K_{i+1} + (I_i \cup K_{n-2i-1})$, where only the edges in the $K_{n-2i-1}$ are allowed to be missing. Let $x$ and $y$ be dominating vertices of $G$. If $G$ had a spanning $x,y$-path then $G-x-y$ would be traceable, but deleting the remaining $i-1$ dominating vertices from $G-x-y$ would yield at least $i+1$ components, so $G-x-y$ is not traceable. Therefore $G$ is not Hamiltonian-connected.
        \item Let $G \in \mathcal{G}^{k}_n(i)$. Then $G$ by definition is a spanning subgraph of $K_{i+k} + (I_i \cup K_{n-2i-k})$, where only the edges in the $K_{n-2i-k}$ are allowed to be missing. Consider any path on $k+1$ dominating vertices and delete its vertices. The result is not traceable because deleting the remaining $i-1$ dominating vertices yields a graph with at least $i+1$ components. Therefore this path of length $k$ is not contained in a Hamilton cycle. Since a path of length $k$ is a linear forest of size $k$, $G$ is not $k$-edge Hamiltonian.
        \item Let $G \in \mathcal{G}^{k}_n(i)$. Then $G$ by definition is a spanning subgraph of $K_{i+k} + (I_i \cup K_{n-2i-k})$, where only the edges in the $K_{n-2i-k}$ are allowed to be missing. Deleting any $k$ of the dominating vertices yields a non-Hamiltonian graph because deleting the remaining $i$ dominating vertices yields a graph with at least $i+1$ components. Therefore $G$ is not $k$-Hamiltonian. \qedhere
    \end{enumerate}
\end{proof}

We will couple \cref{thm:Pdegree} with the following proposition about the graphs $\G{\ell}{i}$.

\begin{prop}\label{prop:uniquedegree}
    Let $n,i\in\Z^+$ and let $\ell\in\Z$. For $-1\le \ell \le n-3$ and $1 \le i \le \frac{n-1-\ell}{2}$, the graph $\G{\ell}{i}$ is the unique graph having nondecreasing degree list 
    \[
        \underbrace{i+\ell, \ldots, i+\ell}_{i \text{ times}}, \underbrace{n-i-1, \ldots, n-i-1}_{n-2i-\ell \text{ times}}, \underbrace{n-1, \ldots, n-1}_{i+\ell \text{ times}}.
    \]
\end{prop}

\begin{proof}
    Let $G$ be a graph having this degree list. Then the $i+\ell$ vertices of degree $n-1$ form a clique $K_{i+\ell}$ and are adjacent to all other vertices of $G$, including those of degree $i+\ell$, which therefore have no other neighbors and form an independent set $I_i$. Each of the remaining $n-2i-\ell$ vertices cannot be adjacent to itself or to the $i$ minimum-degree vertices, so must be adjacent to all $n-i-1$ other vertices in order to have degree $n-i-1$. Therefore $G = K_{i+\ell} + (I_i \cup K_{n-2i-\ell}) = \G{\ell}{i}$. Since $1 \le i \le (n-\ell-1)/2$, we have $i + \ell \le n-i-1 < n-1$.
\end{proof}

\section{Proofs of Theorem \ref{thm:main} and Theorem \ref{thm:extFam}}\label{sec:mainproof}

We are now ready to prove the following.

\begin{thm}\label{thm:gfamily}
     Let $G$ be an $n$-vertex graph with $n\geq 3$ and minimum degree $\delta(G) \ge d$ for some $0 \le d \le (n+\ell-1)/2$ (where $\ell$ is $0, -1, 1, k, k$ in the list below), and let $t\in\{1, \ldots, n-1\}$. 
\begin{enumerate}
\item  If $G$ is not Hamiltonian then $s_t(G) \le \max\set{s_t(\G{0}{i}): \max\set{1,d} \le i \le \tfrac{n-1}{2}}$. 
\item If $G$ is not traceable then $s_t(G) \le \max\set{s_t(\G{-1}{i}): d+1 \le i \le \tfrac{n}{2}}$. 
 \item If $G$ is not Hamiltonian-connected then $s_t(G) \le \max\set{s_t(\G{1}{i}): \max\set{1,d-1} \le i \le \tfrac{n-2}{2}}$. 
 \item If $G$ is not $k$-edge Hamiltonian, for some $k\in\{1, \ldots, n-2\}$, then \\$s_t(G) \le \max\set{s_t(\G{k}{i}): \max\set{1,d-k} \le i \le \tfrac{n-1-k}{2}}$.  
 \item If $G$ is not $k$-Hamiltonian, for some $k\in\{0, \ldots, n-3\}$, then \\$s_t(G) \le \max\set{s_t(\G{k}{i}): \max\set{1,d-k} \le i \le \tfrac{n-1-k}{2}}$. 
 \end{enumerate}
 \end{thm}

\begin{proof}
    Let $G$ be such a graph. Let $\ell$ be the appropriate value from the list $0$, $-1$, $1$, $k$, or $k$, corresponding to the property which $G$ is assumed not to have. Then \cref{thm:Pdegree} implies that there exists an $i^*$ in $\{1,\dots, \floor{(n-1-\ell)/2}\}$ such that 
    \[
        d_j \le 
        \begin{cases}
            i^*+\ell & \text{for } 1 \le j \le i^*\\
            n-i^*-1 & \text{for } i^*+1 \le j \le n-i^*-\ell\\
            n-1 & \text{for } n-i^*-\ell+1 \le j \le n.
        \end{cases}
    \]
    So $d \le \delta(G) = d_1 \le i^* + \ell$. Therefore $i^* \ge d-\ell$.
    
    Let $(c_j)_{j=1}^n$ be the sequence defined by these upper bounds: 
    \[
        c_j = 
        \begin{cases}
            i^*+\ell & \text{for } 1 \le j \le i^*\\
            n-i^*-1 & \text{for } i^*+1 \le j \le n-i^*-\ell\\
            n-1 & \text{for } n-i^*-\ell+1 \le j \le n,
        \end{cases}
    \]
    so $d_j \le c_j$ for every $j$. Notice that, by \cref{prop:uniquedegree}, $c_1 \le \cdots \le c_n$ is the degree list of $\G{\ell}{i^*}$. Therefore, for this value of $i^*$, 
    \[
        s_t(G) = \sum_{j=1}^{n}\binom{d_j}{t} \le \sum_{j=1}^{n}\binom{c_j}{t} = s_t(\G{\ell}{i^*}) \le \max\set{s_t(\G{\ell}{i}): \max\set{1,d-\ell} \le i \le (n-1-\ell)/2}.\qedhere
    \]
\end{proof}

It turns out that the there are only two graphs that we need to consider for the above maximums.

\begin{lem}\label{lem:concavity} For all values of $n \ge 3$, $-1 \le \ell \le n-3$, $1 \le t \le n-1$, and $0\leq d \leq (n+\ell-1)/2$ we have
    \[
        \max\set{s_t(\G{\ell}{i}): \max\set{1,d-\ell} =i_d \le i \le i_0=\floor{(n-1-\ell)/2}} = \max\set{s_t(\GOne{\ell}), s_t(\GTwo{\ell})}.
    \]
\end{lem}

\begin{proof}
    Let $n$, $\ell$, $t$, and $d$ be fixed. For each $i$, let $g_i = s_t(\G{\ell}{i})$ if $t\ge 2$, and let $g_i = 2e(\G{\ell}{i})$ if $t=1$. We show that the sequence $(g_i)_{i=i_d}^{i_0}$ satisfies
    \[
        g_i - g_{i-1} \le g_{i+1} - g_i
    \]
    for every $2 \le i \le \floor{(n-1-\ell)/2}-1$, so the sequence is concave up and maximized at one of the endpoints, $i=i_d$ and $i = i_0$. For each $i$, let $\Delta_i = g_i - g_{i-1}$.
    \begin{align*}
        \Delta_i &= g_i - g_{i-1}\\
        &= i\binom{i+\ell}{t} + (n-2i-\ell)\binom{n-i-1}{t} + (i+\ell)\binom{n-1}{t}\\
        & \quad - \left((i-1)\binom{(i-1)+\ell}{t} + (n-2(i-1)-\ell)\binom{n-(i-1)-1}{t} + ((i-1)+\ell)\binom{n-1}{t}\right)\\
        &= \left(i\binom{i+\ell}{t} - (i-1)\binom{(i-1)+\ell}{t}\right) \\
        & \quad + \left((n-2i-\ell)\binom{n-i-1}{t} - (n-2i+2-\ell)\binom{n-i}{t}\right)\\
        & \quad + \left((i+\ell)\binom{n-1}{t} - ((i-1)+\ell)\binom{n-1}{t}\right)\\
        &= \left(\binom{i+\ell}{t} + (i-1)\binom{(i-1)+\ell}{t-1}\right) + \left(-(n-2i-\ell)\binom{n-i-1}{t-1} - 2\binom{n-i}{t}\right)+ \binom{n-1}{t},
    \end{align*}
    where the last step follows from using Pascal's identity in two places. As $\binom{x}{t}$ is a weakly increasing function of $x$ for all $x$, by comparing term by term we have
    \begin{align}
        &\Delta_i - \binom{n-1}{t}= \left(\binom{i+\ell}{t} + (i-1)\binom{(i-1)+\ell}{t-1}\right) + \left(-(n-2i-\ell)\binom{n-i-1}{t-1} - 2\binom{n-i}{t}\right)\nonumber\\
        &\le \left(\binom{(i+1)+\ell}{t} + i\binom{i+\ell}{t-1}\right) + \left(-(n-2(i+1)-\ell)\binom{n-(i+1)-1}{t-1} - 2\binom{n-(i+1)}{t}\right)\label{eq:delta}\\
        &= \Delta_{i+1} - \binom{n-1}{t}\nonumber.\qedhere
    \end{align}

    \cref{thm:main} follows from \cref{thm:gfamily} and \cref{lem:concavity}. Now we turn our attention to determining the set of all extremal graphs.

    \begin{lem}\label{lem:maximizers}
        In the same setting as \cref{lem:concavity}, if $i \notin \set{i_d,i_0}$ then \[s_t(\G{\ell}{i}) < \max\set{s_t(\GOne{\ell}), s_t(\GTwo{\ell})}.\]
    \end{lem}
    \begin{proof}
        With hypotheses as in \cref{lem:concavity}, we first show the following claim: for $i \le n-t$ we have $\Delta_i < \Delta_{i+1}$, and otherwise $\Delta_i = \Delta_{i+1} = \binom{n-1}{t} \ge 1$. To prove this claim, notice that if $i \le n-t$ then $\binom{n-i}{t}-\binom{n-i-1}{t} = \binom{n-i-1}{t-1} \ge 1$, so the inequality $\Delta_i \le \Delta_{i+1}$ is strict by using the fourth term in \cref{eq:delta}. Otherwise using $i \le (n-1-\ell)/2$ and $i \ge n-t+1$ we have
        \begin{align*}
            2i + \ell + 1 \le n &\le i+t-1\\
            i + \ell + 1 &\le t-1
        \end{align*}
        so all four terms of $\Delta_i-\binom{n-1}{t}$ and $\Delta_{i+1}-\binom{n-1}{t}$ are zero. This concludes the proof of the claim.
    
        Now we prove the lemma by contradiction. Suppose that there is some $i \notin \set{i_d,i_0}$ such that $g_i = s_t(\G{\ell}{i}) = \max\set{s_t(\GOne{\ell}), s_t(\GTwo{\ell})} = \max\set{g_{i_d},g_{i_0}}$. We show that then $g_{i_d} = g_{i_d+1} = \cdots = g_{i_0-1} = g_{i_0}$. 
        
        First, consider the values $i_d, i, i_0$. If $g_{i_d} = g_i$ is the maximum, then by weak convexity
        \[
            0 = \frac{g_{i} - g_{i_d}}{i-i_d} \le \frac{g_{i_0} - g_i}{i_0 - i},
        \]
        so $g_{i_0}\ge g_i$, but $g_i$ is the maximum, so $g_{i_0} = g_i = g_{i_d}$. A symmetric argument shows that if $g_{i_0} = g_i$ is the maximum, then $g_{i_d} = g_i = g_{i_0}$. Therefore we assume going forward that $g_{i_d} = g_i = g_{i_0}$ is the maximum.
        
        For all other values $j \notin \set{i_d, i, i_0}$, either $i_d < j < i$ or $i < j < i_0$. We address the first case, and the second is similar. When $i_d < j < i$, by weak concavity we have
        \[
            \frac{g_i - g_j}{i-j} \le \frac{g_{i_0} - g_i}{i_0 - i} = 0,
        \]
        so $g_i \le g_j$, and $g_i$ is the maximum, so $g_i = g_j$. Therefore $g_{i_d} = g_{i_d+1} = \cdots = g_{i_0-1} = g_{i_0}$.

        Then $\Delta_{i_{d}+1} = \Delta_{i_d+2} = \cdots = \Delta_{i_0-1} = \Delta_{i_0} = 0$, contradicting the fact that two consecutive values $\Delta_k$ and $\Delta_{k+1}$ cannot both equal $0$ (from the first claim of this proof). From the contradiction we conclude that there is no $i \notin \set{i_d,i_0}$ such that $s_t(\G{\ell}{i}) = \max\set{s_t(\GOne{\ell}), s_t(\GTwo{\ell})}$.
    \end{proof}  
\end{proof}

Now we can prove \cref{thm:extFam}.

\extremalgraphs*

\begin{proof}
    Let $G$ be an $n$-vertex graph with $\delta(G)\geq d$ which maximizes $s_t(G)$ subject to not being one of Hamiltonian, traceable, Hamiltonian-connected, $k$-edge Hamiltonian, or $k$-Hamiltonian. Let $\ell$ equal $0, -1, 1, k$, or $k$ depending on which of these properties we are considering, respectively. 
    
    Let $d_1 \leq d_2 \leq \cdots \leq d_n$ be the degree sequence of $G$. As in the proof of \cref{thm:gfamily}, by \cref{thm:Pdegree} there exists an $i^*$, $i_d \le i^* \le i_0$, such that 
    \[
        d_j \le 
        \begin{cases}
            i^*+\ell & \text{for } 1 \le j \le i^*\\
            n-i^*-1 & \text{for } i^*+1 \le j \le n-i^*-\ell\\
            n-1 & \text{for } n-i^*-\ell+1 \le j \le n.
        \end{cases}
    \]
    For $1\leq j\leq n$ define
    \[
        c_j = 
        \begin{cases}
            i^*+\ell & \text{for } 1 \le j \le i^*\\
            n-i^*-1 & \text{for } i^*+1 \le j \le n-i^*-\ell\\
            n-1 & \text{for } n-i^*-\ell+1 \le j \le n.
        \end{cases}
    \]
Then $d_j\leq c_j$ for all $1\leq j\leq n$. Moreover, by the argument given in the proof of \cref{thm:gfamily}, coupled with the result of \cref{lem:concavity}, we know that
 \begin{equation}\label{binomSum}
        s_t(G) =\begin{cases} \displaystyle\sum_{j=1}^{n}\binom{d_j}{t} = \displaystyle\sum_{j=1}^{n}\binom{c_j}{t} & \text {for } t\ge 2 \\
        \displaystyle\frac{1}{2}\sum_{j=1}^{n}d_j = \frac{1}{2}\sum_{j=1}^{n}c_j & \text{for } t=1
        \end{cases} \quad = \max\set{s_t(\GOne{\ell}), s_t(\GTwo{\ell})}.
    \end{equation}
Since Proposition \ref{prop:properties} assures us that $\GOne{\ell}$ and $\GTwo{\ell}$ do not have their associated property (Hamiltonian for $\ell=0$, traceable for $\ell=-1$, Hamiltonian-connected for $\ell=1$, $k$-edge Hamiltonian or $k$-Hamiltonian for $\ell=k$), we know that $\GOne{ \ell}$ or $\GTwo{\ell}$ (or both) is a member of the extremal family we are looking for. By  
\cref{prop:uniquedegree} and \cref{lem:maximizers}, they are the only possible members with degree sequence $c_1 \le \cdots \le c_n$.

Suppose $G$ is not equal to either of $\GOne{\ell}$ and $\GTwo{\ell}$. Let $v_1, \ldots, v_n$ be the vertices of $G$, with $d(v_j)=d_j$ for all $j$. Let $w_1, \ldots, w_n$ be the vertices of $\GOne{\ell}$ or $\GTwo{\ell}$ (whichever maximizes $s_t(G)$), with $d(w_j)=c_j$ for all $j$. Since $d_j\leq c_j$ for all $j$, \cref{binomSum} tells us that $d_j=c_j$ whenever $c_j\geq t$, that is, whenever the binomial coefficient in the sum is nonzero. Recall that $\GOne{\ell}$ and $\GTwo{\ell}$ are defined to be $K_{i+\ell}+(I_i\cup K_{n-2i-\ell})$, where $i=\max\{1,d-\ell\}$ and $i=\lfloor\tfrac{n-1-\ell}{2}\rfloor$, respectively. For convenience, we let  $A, I, B$ denote the vertices in each of the three parts of this graph, namely the $K_{i+\ell},I_i$, and $K_{n-2i-\ell}$, respectively. Note that each vertex in $A=K_{i + \ell}$ is adjacent to all others in the graph, and so $n-1\geq t$ implies that $d(v_j)=d(w_j)=n-1$ for all $j$ such that $w_j\in A$. 

The vertices $v_j$ corresponding to $w_j\in A$ force every other vertex in $G$ to have degree at least $i+\ell$. But we know that in $\GOne{\ell}, \GTwo{\ell}$, there are $i$ vertices whose degree is equal to $i+\ell$. Since $d_j\leq c_j$ for all $j$, this means there exist $i$ vertices in $G$ whose degrees are exactly equal to $i+\ell$ as well.

There are $n-2i-\ell$ vertices of $G$ whose degrees are yet to be determined. These correspond to the vertices $w\in B$, and we know that for each such $w$, $d(w)=n-i-1$. 
Our remaining vertices in $G$ can certainly have degrees no larger than this. If $n-i-1\geq t$, then the remaining vertices in $G$ must have all degrees exactly equal to $n-i-1$, due to (\ref{binomSum}). However, if $n-i-1 <t$, then while the remaining vertices must all be adjacent to those $v_j$ corresponding to the $A$-vertices, and cannot be adjacent to any of those $v_j$ corresponding to the $I$-vertices, their adjacencies among themselves can be anything and they will still optimize $s_t(G)$.
\end{proof}

\begin{example}\label{ex:otherextremal}
    Let $n=10$, $\ell=0$, and $d=4$. Notice that $d \le (n+\ell-1)/2$, so this choice of $d$ is valid. Then $i_d = 4$ and $i_0=4$, so the extremal graphs are all in $\mathcal{G}^0_{10}(4)$ by \cref{thm:extFam}. For every $6 \le t \le 9$, we have $n-i_d-1 = n-i_0-1 = 5 < t$, so there are multiple extremal graphs: $K_4+(I_4\cup I_2)$ and $K_4+(I_4\cup K_2)$ have the same numbers of $t$-stars.
\end{example}

\section{Connectedness}\label{sec:kconn}

The property of $k$-connectedness is $(n+k-2)$-stable \cite[Theorem 9.7]{BC76}, which yields a Chv\'{a}tal-like degree condition by \cref{thm:Pdegree}. However, this condition is not best possible for $k$-connectedness because, for $i > 1$, the graph $\G{k-2}{i}$ has at least $k$ dominating vertices so is $k$-connected. Therefore, we address $k$-connectedness separately to obtain a tight upper bound using different extremal graphs.

Recall from the introduction that for $n, k, i \in \mathbb{Z}$, where $1 \le k \le n-2$ and $1 \le i \le \frac{n-k+1}{2}$, $\h{k}{i}:=K_{k-1} + (K_i \cup K_{n-k-i+1})$, and $\h{k}{i}$ is not $k$-connected. The graph $\h{k}{i}$ has degree list \[\underbrace{i+k-2, \ldots, i+k-2}_{i \text{ times}}, \underbrace{n-i-1, \ldots, n-i-1}_{n-k-i+1 \text{ times}}, \underbrace{n-1, \ldots, n-1}_{k-1 \text{ times}}.\]

We use the following theorem, which is stronger than the Chv\'atal-like degree condition for $k$-connectedness guaranteed by \cref{thm:Pdegree}.
\begin{thm}[Bondy \cite{Bondy69}, Boesch \cite{Boesch74}]\label{thm:connecteddegree}
    Let $G$ be an $n$-vertex graph and $d_1 \le \cdots \le d_n$ its degrees. For $1 \le k \le n-2$, if $G$ is not $k$-connected, then there is an integer $1 \le i \le \frac{n-k+1}{2}$ for which $d_i \le i+k-2$ and $d_{n-k+1} \le n-i-1$. In other words, $G$ has at least $i$ vertices of degree at most $i+k-2$ and at most $k-1$ vertices of degree at least $n-i$.
\end{thm}

First we show that there is an extremal graph in the $\h{k}{i}$ family.

\begin{thm}\label{thm:hfamily}
    Let $G$ be an $n$-vertex graph with $n \ge 3$ and minimum degree $\delta(G) \ge d$ for some $0 \le d \le (n+k-3)/2$, where $1 \le k \le n-2$, and let $t \in \set{1, \ldots, n-1}$. If $G$ is not $k$-connected, then $s_t(G) \le \max\set{s_t(\h{k}{i}): \max\set{1,d-k+2} \le i \le \tfrac{n-k+1}{2}}$.
\end{thm}

\begin{proof}
    Let $G$ be such a graph. Then \cref{thm:connecteddegree} implies that there exists an $i^*$ in $1 \le i^* \le (n-k+1)/2$ such that 
    \[
        d_j \le 
        \begin{cases}
            i^*+k-2 & \text{for } 1 \le j \le i^*\\
            n-i^*-1 & \text{for } i^*+1 \le j \le n-k+1\\
            n-1 & \text{for } n-k+2 \le j \le n.
        \end{cases}
    \]
    So $d \le \delta(G) = d_1 \le i^*+k-2$. Therefore $i^* \ge d-k+2$, and $i^* \ge \max\set{1,d-k+2} =: i_d$. Let $(c_j)_{j=1}^n$ be the sequence defined by these upper bounds: 
    \[
        c_j = 
        \begin{cases}
            i^*+k-2 & \text{for } 1 \le j \le i^*\\
            n-i^*-1 & \text{for } i^*+1 \le j \le n-k+1\\
            n-1 & \text{for } n-k+2 \le j \le n,
        \end{cases}
    \]
    so $d_j \le c_j$ for every $j$, and $c_1 \le \cdots \le c_n$ is the degree list of $\h{k}{i^*}$. Therefore, we have either
    \[
        s_t(G) = \sum_{j=1}^{n}\binom{d_j}{t} \le \sum_{j=1}^{n}\binom{c_j}{t} = s_t(\h{k}{i^*}) \le \max\set{s_t(\h{k}{i}): i_d \le i \le (n-k+1)/2}
    \]
    for $t \ge 2$, or, similarly,
    \[
        e(G) = s_1(G) = 2\sum_{j=1}^{n}d_j \le 2\sum_{j=1}^{n}c_j = s_1(\h{k}{i^*}) \le \max\set{s_1(\h{k}{i}): i_d \le i \le (n-k+1)/2}
    \]
    for $t=1$.
\end{proof}

Now we find the maximum number of $t$-stars within the $\h{k}{i}$ family.

\begin{lem}\label{lem:connectedconcavity} For all values of $n \ge 3$, $1\le k \le n-2$, $0 \le d \le (n+k-3)/2$, and $1 \le t \le n-1$, let $i_d = \max\set{1,d-k+2}$ and $i_0 = \floor{(n-k+1)/2}$. Then 
    \[
        \max\set{s_t(\h{k}{i}): i_d \le i \le i_0} =
        \begin{cases}
            s_t(\h{k}{i_d}) & \text{for } t=1\\
            \max\set{s_t(\h{k}{i_d}), s_t(\h{k}{i_0})} & \text{for }2 \le t \le n-1.
        \end{cases}
    \]
\end{lem}

\begin{proof}
    Let $n$, $k$, and $t$ be fixed. For each $i$ in $i_d \le i \le i_0$, let $h_i = s_t(\h{k}{i})$ if $t \ge 2$, and let $h_i = 2e(\h{k}{i})$ if $t=1$. We show that the sequence $(h_i)_{i=i_d}^{i_0}$ satisfies
    \[
        h_i - h_{i-1} \le h_{i+1} - h_i
    \]
    for every $i_d+1 \le i \le i_0-1$, so the sequence is concave up and maximized at one of the endpoints, $i=i_d$ and $i = i_0$. For ease of notation, for each $i$, let $\Delta_i = h_i - h_{i-1}$. For $t=1$ we also show that $\Delta_i < 0$, so the sequence $(h_i)_{i=i_d}^{i_0}$ is decreasing and maximized at the left endpoint $i=i_d$.
    \begin{align}
        \Delta_i &= h_i - h_{i-1}\nonumber\\
        &= i\binom{i+k-2}{t} + (n-k-i+1)\binom{n-i-1}{t} + (k-1)\binom{n-1}{t}\nonumber\\
        & \quad - \left((i-1)\binom{(i-1)+k-2}{t} + (n-k-(i-1)+1)\binom{n-(i-1)-1}{t} + (k-1)\binom{n-1}{t}\right)\nonumber\\
        &= \left(i\binom{i+k-2}{t} - (i-1)\binom{i+k-3}{t}\right) \nonumber\\
        & \quad + \left((n-k-i+1)\binom{n-i-1}{t} - (n-k-i+2)\binom{n-i}{t}\right)\nonumber\\
        & \quad + \left((k-1)\binom{n-1}{t} - (k-1)\binom{n-1}{t}\right)\nonumber\\
        &= \left(\binom{i+k-2}{t} + (i-1)\binom{i+k-3}{t-1}\right)  + \left(-(n-k-i+1)\binom{n-i-1}{t-1} - \binom{n-i}{t}\right)\label{eq:hdelta},
    \end{align}
    where the last step follows from using Pascal's identity in two places.
    
    When $t=1$ \cref{eq:hdelta} simplifies to 
    \[
        \Delta_i = (i+k-2) + (i-1) - (n-k-i+1) - (n-i) = 4i+2k-2n-4,
    \]
    and using the fact that $i \le (n-k+1)/2$, we have $4i+2k-2n \le 2$, so $\Delta_i \le -2 < 0$.
    
    For $1 \le t \le n-1$, as $\binom{x}{t}$ is a weakly increasing function of $x$ for all $x$, we have
    \begin{align*}
        \Delta_i &= \left(\binom{i+k-2}{t} + (i-1)\binom{i+k-3}{t-1}\right) + \left(-(n-k-i+1)\binom{n-i-1}{t-1} - \binom{n-i}{t}\right)\\
        &\le \left(\binom{(i+1)+k-2}{t} + i\binom{i+k-2}{t-1}\right) + \left(-(n-k-i)\binom{n-i-2}{t-1} - \binom{n-i-1}{t}\right)\\
        &= \Delta_{i+1}.\qedhere\\
    \end{align*}
\end{proof}

\cref{thm:kconn} follows from \cref{thm:hfamily} and \cref{lem:connectedconcavity}. Here we use the fact that the graph $\h{k}{i_d}$ or the graph $\h{k}{i_0}$ achieves the upper bound; both graphs are not $k$-connected, as proved above \cref{thm:kconn} in the introduction.

\section{Proof Of Theorem \ref{l0conj}}\label{sec:bigConj}

In this section, we work to identify which of $\G{\ell}{1}$ and $\GTwo{\ell}$ contains more $t$-stars, depending on the value of $t$ with respect to $n$ and $\ell$.  Our results constitute a proof of Theorem \ref{l0conj}.

\subsection{Large \textit{t}}\label{subsec:larget}

First, we prove that $\GTwo{\ell}$ contains the maximum number of $t$-stars when $t\geq (n+\ell+1)/2$.

\begin{prop}\label{prop:max} Let $n,i,\ell\in\mathbb{Z}$, $-1\leq \ell\leq n-3$, and $\tfrac{n+\ell+1}{2} \leq t \leq n-1$. Then $$s_t(\G{\ell}{1}) \leq s_t(\GTwo{\ell}).$$ 
Moreover, when $0 \le \ell \le n-5$, the inequality is strict.
\end{prop}

\begin{proof} 
By Lemma \ref{lem:concavity}, the sequence of $s_t(\G{\ell}{i})$ values for $1\leq i \leq i_0 = \floor{\frac{n-\ell-1}{2}}$ is concave up, which means it is maximized at one of the endpoints. 
The definitions of $\G{\ell}{1}$ and $\GTwo{\ell}$ tell us that, when $t\geq 2$,
\begin{equation}\label{G1} 
s_t(\G{\ell}{1})= \binom{\ell+1}{t}+ (\ell+1)\binom{n-1}{t}+ (n-\ell-2)\binom{n-2}{t}
\end{equation}
and
\begin{equation}\label{G2}   
s_t(\GTwo{\ell})=i_0\binom{i_0 +\ell}{t} + \left(i_0+\ell\right)\binom{n-1}{t} + \left(n-2i_0-\ell\right)\binom{n-i_0 -1}{t}.
\end{equation}
When $t=1$, these values are multiplied by $1/2$, and otherwise the same argument holds. We can simplify these expressions significantly using the following two observations. First, $t>\ell+1$ because $\ell+1< t \Leftarrow \ell+1 <\tfrac{n+\ell+1}{2} \Leftrightarrow \ell+1 < n.$ Second, $t> \lceil \tfrac{n-\ell-1}{2}\rceil +\ell$ because $ t> \lceil \tfrac{n-\ell-1}{2}\rceil +\ell \Leftarrow  t>\tfrac{n-\ell}{2}+\ell \Leftrightarrow t>\tfrac{n+\ell}{2}$. 
Thus we can disregard the first term of (\ref{G1}) and the first term of (\ref{G2}). In fact, we can also disregard the third term of (\ref{G2}) since 
$n-\lfloor\tfrac{n-\ell-1}{2}\rfloor -1=\lceil\tfrac{n-\ell-1}{2}\rceil+ \ell$. So we get that 
$$s_t(\GTwo{\ell}) - s_t(\G{\ell}{1})= \left(\left\lfloor \tfrac{n-\ell-1}{2}\right\rfloor-1\right)\binom{n-1}{t}- (n-\ell-2)\binom{n-2}{t}.$$

Notice that when $t=n-1$, the difference above is $\left\lfloor \tfrac{n-\ell-1}{2}\right\rfloor-1 \geq 0$ as desired, and the inequality is strict unless $\ell\ge n-4$. Thus, for the remainder of the proof we may assume that $t\leq n-2$. 

Using the fact that $t < n-1$, we can combine the two binomial terms into one as follows:
\begin{eqnarray*} s_t(\GTwo{\ell}) - s_t(\G{\ell}{1}) 
&=& \left( \lfloor\tfrac{n-\ell-1}{2}\rfloor -1 \right) \binom{n-2}{t} \frac{(n-1)}{(n-1-t)}-(n-\ell-2)\binom{n-2}{t}\\
&=& \binom{n-2}{t}\left(( \lfloor\tfrac{n-\ell-1}{2}\rfloor -1 )\frac{(n-1)}{(n-1-t)}-(n-\ell-2)\right).
\end{eqnarray*}

Since $t$ is an integer, we have $t \ge \ceil{\frac{n+\ell+1}{2}}$. Using the fact that $\floor{\frac{n-\ell-1}{2}} + \ceil{\frac{n+\ell+1}{2}} = n$,
\begin{eqnarray*} s_t(\GTwo{\ell}) - s_t(\G{\ell}{1}) &\geq& 
 \binom{n-2}{t}\left(\left( \floor[\bigg]{\frac{n-\ell-1}{2}} -1 \right)\tfrac{(n-1)}{(n-1-\ceil{(n+\ell+1)/2})}-(n-\ell-2)\right)\\
&=& \binom{n-2}{t}\left((n-1)-(n-\ell-2)\right)=\binom{n-2}{t}\left(\ell+1\right) \ge 0,
\end{eqnarray*}
since $\ell \ge -1$, with a strict inequality for all $\ell \ge 0$.
\end{proof}

\subsection{Small \textit{t}}\label{subsec:smallt}

Next, we use an inductive argument to show that in the non-Hamiltonian case (i.e. when $\ell=0$), $\G{0}{1}$ contains the maximum number of $t$-stars when $t< (n+1)/2$. 

\begin{prop}\label{prop:smallt} Let $n\in\mathbb{Z}$, $n\geq 4$, $\ell=0$, and $1\leq t\leq n/2$. Then $$s_t(\G{0}{1}) \geq s_t(\GTwo{0}).$$ Moreover, when $n \ge 6$, the inequality is strict.
\end{prop}

\begin{proof}
First suppose $t=1$. The graph $\G{0}{1}$ is $K_{n+1}$ with a pendent edge so has $e(\G{0}{1}) = \binom{n-1}{2}+1 = (n^2 - 3n + 4)/2$.

    The graph $\GTwo{0}$ when $n$ is odd is $K_{(n-1)/2} + I_{(n+1)/2}$ so has $e(\GTwo{0}) = \binom{(n-1)/2}{2} + \frac{n-1}{2}\cdot \frac{n+1}{2} =(3n^2-4n+1)/8$. In this case $e(\G{0}{1}) - e(\GTwo{0}) = (n^2 - 8 n + 15)/8$, which is nonnegative for odd $n \ge 3$ and positive for all $n \ge 6$.

    The graph $\GTwo{0}$ when $n$ is even is $K_{n/2-1} + (I_{n/2-1}\cup K_2)$ so $e(\GTwo{0}) = (n/2-1)(n/2-1) + \binom{n/2-1}{2} + 2(n/2-1)+1 = (3n^2 - 6n + 8)/8$. In this case $e(\G{0}{1}) - e(\GTwo{0}) =  (n^2 - 6 n + 8)/8$, which is nonnegative for even $n \ge 2$ and positive for all $n \ge 5$.
    
    Notice that the weak inequality for $t=1$ can alternatively be obtained by \cref{thm:edgemax} because the graph $\GTwo{0}$ is not Hamiltonian (an $n$-stable property) by \cref{prop:properties} and has $n \ge 3$.

Now we suppose $t\ge 2$. Notice that
$$s_t(\G{0}{1})=(n-2)\binom{n-2}{t}+ \binom{n-1}{t}.$$
When $n$ is odd,
$$s_t(\GTwo{0})=s_t(\G{0}{\tfrac{n-1}{2}})=\frac{n+1}{2}\binom{(n-1)/2}{t}+ \frac{n-1}{2}\binom{n-1}{t},$$
and when $n$ is even,
$$s_t(\GTwo{0})=s_t(\G{0}{\tfrac{n}{2}-1})=\left(\frac{n}{2}-1\right)\binom{n/2-1}{t}+ \left(\frac{n}{2}-1\right)\binom{n-1}{t}+2\binom{n/2}{t}.$$

We will first take care of the $t=2$ case. Using the equations above, $s_2(\G{0}{1})=(n-2)(n^2-4n+5)/2$, $s_2(\GTwo{0})=(n-1)(5n^2-14n+5)/16$ when $n$ is odd, and $s_2(\GTwo{0})=(n-2)(5n^2-14n+16)/16$ when $n$ is even. It is easy to see that $s_2(\G{0}{1})> s_2(\GTwo{0})$ for all $n\geq 6$.

Now suppose $t\ge 3$. Define a function $f$ that measures the gap between $s_t(\G{0}{1})$ and $s_t(\GTwo{0})$:
 \begin{align*}
    f(n,t)&= s_t(\G{0}{1}) - s_t(\GTwo{0}) \\
    &=
    \begin{cases}
      s_t(\G{0}{1}) - s_t(\G{0}{\tfrac{n-1}{2}}), & \text{if}\ n \text{ is odd} \\
      s_t(\G{0}{1}) - s_t(\G{0}{\tfrac{n}{2}-1}), & \text{if}\ n \text{ is even}\end{cases}
  \end{align*}

We can simplify $f(n,t)$ in the following way. When $n$ is odd,
 \begin{align}
    f(n,t)&= s_t(\G{0}{1}) - s_t(\G{0}{\tfrac{n-1}{2}}) \nonumber\\
    &= \left[(n-2){{n-2} \choose t}+ {{n-1} \choose t}\right] - \left[\frac{n+1}{2}{{(n-1)/2} \choose t}+ \frac{n-1}{2}{{n-1}\choose t}\right]\nonumber\\
    &=\left[\frac{(n-2)(n-t-1)}{n-1}{{n-1} \choose t}+ {{n-1} \choose t}\right] - \left[\frac{n+1}{2}{{(n-1)/2} \choose t}+ \frac{n-1}{2}{{n-1}\choose t}\right]\nonumber\\
    &={{n-1} \choose t}\left[\frac{(n-2)(n-t-1)}{n-1}+1-\frac{n-1}{2}\right]-{{(n-1)/2} \choose t}\frac{n+1}{2}\nonumber\\
     &={{n-1} \choose t}\left[\frac{(n-1)}{2}-t+\frac{t}{n-1}\right]-{{(n-1)/2} \choose t}\frac{n+1}{2},\label{eq:odd}
    \end{align}
and when $n$ is even,
 \begin{align}
    f(n,t)&= s_t(\G{0}{1}) - s_t(\G{0}{\tfrac{n}{2}-1}) \nonumber\\
    &= \left[(n-2){{n-2} \choose t}+ {{n-1} \choose t}\right] - \left[\left(\frac{n}{2}-1\right){{n/2-1} \choose t}+ \left(\frac{n}{2}-1\right){{n-1} \choose t}+2{{n/2} \choose t}\right]\nonumber\\
    &={{n-1} \choose t}\left[\frac{(n-2)(n-t-1)}{n-1}+1-\left(\frac{n}{2}-1\right)\right] - \left[\left(\frac{n}{2}-1\right){{n/2-1} \choose t}+2{{n/2} \choose t}\right]\nonumber\\
    &={{n-1} \choose t}\left[\frac{n}{2}-t+\frac{t}{n-1}\right] - {{n/2} \choose t}\left[\left(\frac{n}{2}-1\right)\frac{n/2-t}{n/2}+2\right]\nonumber\\
        &={{n-1} \choose t}\left[\frac{n}{2}-t+\frac{t}{n-1}\right] - {{n/2} \choose t}\left[ \frac{n}{2}-t+1+\frac{2t}{n}\right].\label{eq:even}
    \end{align}
    
We will now use an inductive argument to show that $s_t(\G{0}{1}) \geq s_t(\GTwo{0})$ for all pairs of values $n,t\in\mathbb{Z}$ where $n\geq 6$ and $t\in\{3,\dots,\floor{n/2}\}$. For the base case we prove Claim \ref{claim:l=0base}, which states that $f(2t,t)$ is 
positive for all $t\geq 3$.  For the inductive step, we prove Claim \ref{claim:l=0induct}, which states that for fixed $t\geq 2$ and $n\geq 2t$, $f(n,t)$ is strictly increasing with respect to $n$.

\begin{claim}\label{claim:l=0base}
    Let $t\in\mathbb{Z}$, $t\geq 3$. Then $f(2t,t)>0$.
\end{claim}
\begin{proofc} 
    Note that $2t$ is even, so by \cref{eq:even},
    \begin{align*}
        f(2t,t) &= \binom{2t-1}{t}\left[\frac{2t}{2}-t+\frac{t}{2t-1}\right]-\binom{\frac{2t}{2}}{t}\left[\frac{2t}{2}-t+1+\frac{2t}{2t}\right] \\
        &= \binom{2t-1}{t}\cdot\frac{t}{2t-1}-2= \binom{2t-2}{t-1}-2> \binom{2t-2}{1}-2 = 2t-4 >0.\qedhere
    \end{align*}
\end{proofc}

\begin{claim}\label{claim:l=0induct}
    Let $n,t\in\mathbb{Z}$, $t\geq 3$, and $n\geq 2t$.  Then $f(n+1,t)-f(n,t)>0$. 
\end{claim}
\begin{proofc} 
We will consider two cases based on the parity of $n$.

\begin{case}
Suppose $n$ is even. Using Equations (\ref{eq:odd}) and (\ref{eq:even}), 
\begin{align*}
    &f(n+1,t)-f(n,t)\\
    &=\left[{{n} \choose t}\left[\frac{n}{2}-t+\frac{t}{n}\right]-{{n/2} \choose t}\frac{n+2}{2}\right]-\left[{{n-1} \choose t}\left[\frac{n}{2}-t+\frac{t}{n-1}\right] - {{n/2} \choose t}\left[ \frac{n}{2}-t+1+\frac{2t}{n}\right]\right] \\
    &={{n-1} \choose t}\left[\frac{n}{n-t}\left(\frac{n}{2}-t+\frac{t}{n}\right)-\left(\frac{n}{2}-t+\frac{t}{n-1}\right)\right]-{{n/2} \choose t}\left[\frac{n+2}{2}-\left(\frac{n}{2}-t+1+\frac{2t}{n}\right)\right]\\
    &={{n-1} \choose t}t\left[\frac{n/2-t+1}{n-t}-\frac{1}{n-1}\right]-{{n/2} \choose t}t\left(1-\frac{2}{n}\right)\\
    &>{{n-1} \choose t}t\left[\frac{n/2-t+1}{n-t}-\frac{1}{n-1}\right]-{{n/2} \choose t}t\\
    &=t\left[{{n-1} \choose t}\left(\frac{n/2-t+1}{n-t}-\frac{1}{n-1}\right)-{{n/2} \choose t}\right]  
  \end{align*}
We will show that ${{n-1} \choose t}\left(\frac{n/2-t+1}{n-t}-\frac{1}{n-1}\right)-{{n/2} \choose t}\ge 0$ by proving the equivalent statement
\begin{equation*}
   \frac{ {{n-1} \choose t}}{{{n/2} \choose t}}\left(\frac{n/2-t+1}{n-t}-\frac{1}{n-1}\right)\ge 1.
\end{equation*}
Notice that
\begin{align}
    \frac{ {{n-1} \choose t}}{{{n/2} \choose t}}&\left(\frac{n/2-t+1}{n-t}-\frac{1}{n-1}\right)\nonumber\\
    &= \frac{(n-1)!}{(n-t-1)!}\frac{(n/2-t)!}{(n/2)!}\left(\frac{(n/2-t+1)(n-1)-n+t}{(n-t)(n-1)}\right)\nonumber\\
    &=\frac{n-t}{n/2-t+1}\cdot\frac{n-t+1}{n/2-t+2}\cdots\frac{n-3}{n/2-2}\cdot\frac{n-2}{n/2-1}\cdot\frac{n-1}{n/2}\cdot\left(\frac{(n/2-t+1)(n-1)-n+t}{(n-t)(n-1)}\right)\nonumber\\
     &=\frac{\cancel{n-t}}{n/2-t+1}\cdot\left[\frac{n-t+1}{n/2-t+2}\cdots\frac{n-3}{n/2-2}\right]\cdot\frac{n-2}{n/2-1}\cdot\frac{\cancel{n-1}}{n/2}\cdot\left(\frac{(n/2-t+1)(n-1)-n+t}{\cancel{(n-t)}\cancel{(n-1)}}\right)\nonumber\\
     &=\left[\frac{n-t+1}{n/2-t+2}\cdots\frac{n-3}{n/2-2}\right]\cdot\frac{2\cancel{(n-2)}}{\cancel{n-2}}\cdot\left(\frac{(n/2-t+1)(n-1)-n+t}{(n/2-t+1)(n/2)}\right)\label{eq:1}
\end{align}
Notice that we are using the fact that $t \ge 3$ to ensure that the $n-t$ and $n-2$ factors above are distinct. Further, 
\begin{equation}\label{eqn4}
\frac{n-t+i}{n/2-t+i+1}=\frac{n/2-t+i+1+n/2-1}{n/2-t+i+1}=1+\frac{n/2-1}{n/2-t+i+1}>1
\end{equation}
for any $i\geq 0$ since $t\leq n/2$ and $n\geq 4$. Thus, the term in brackets of the last line in \cref{eq:1} above, by \cref{eqn4}, is strictly greater than 1. Thus,
\begin{align*}
    \frac{ {{n-1} \choose t}}{{{n/2} \choose t}}&\left(\frac{n/2-t+1}{n-t}-\frac{1}{n-1}\right)\\
     &=\left[\frac{n-t+1}{n/2-t+2}\cdots\frac{n-3}{n/2-2}\right]\cdot2\cdot\left(\frac{(n/2-t+1)(n-1)-n+t}{(n/2-t+1)(n/2)}\right)\\
     &>2\cdot\frac{(n/2-t+1)(n-1)-n+t}{(n/2-t+1)(n/2)}\\
     &=2\left(\frac{n-2}{n/2}-\frac{n/2-1}{(n/2-t+1)(n/2)}\right)\\
     &\geq 2\left(\frac{n-2}{n/2}-1+\frac{2}{n}\right)=2\left(1-\frac{2}{n}\right) \ge 1.
\end{align*}

This last line follows from the fact that $t\leq \frac{n}{2}$, and the last inequality holds for $n\geq 4$. This completes the proof for the case when $n$ is even. 
\end{case}

\begin{case}
 Suppose $n$ is odd. Using Equations (\ref{eq:odd}) and (\ref{eq:even}),  
\begin{align*}
    f(n+1,t)-f(n,t)
    &=\left[{{n} \choose t}\left[\frac{n+1}{2}-t+\frac{t}{n}\right] - {{(n+1)/2} \choose t}\left[ \frac{n+1}{2}-t+1+\frac{2t}{n+1}\right]\right]\\
    &\qquad\qquad\qquad -\left[{{n-1} \choose t}\left[\frac{n-1}{2}-t+\frac{t}{n-1}\right]-{{(n-1)/2} \choose t}\frac{n+1}{2}\right]\\
    &={{n-1} \choose t}\left[\frac{n}{n-t}\left(\frac{n+1}{2}-t+\frac{t}{n}\right)-\left(\frac{n-1}{2}-t+\frac{t}{n-1}\right)\right]\\
    &\qquad\qquad\qquad -{{(n+1)/2} \choose t}\left[\frac{n+1}{2}-t+1+\frac{2t}{n+1}-\frac{(n+1)/2-t}{(n+1)/2}\cdot\frac{n+1}{2}\right]\\
    &={{n-1} \choose t}\left[\frac{(n-2t+3)t}{2(n-t)}+\frac{n-t-1}{n-1}\right]
    -{{(n+1)/2} \choose t}\left(1+\frac{2t}{n+1}\right)\\
    &>{{n-1} \choose t}\frac{(n-2t+3)t}{2(n-t)}
    -{{(n+1)/2} \choose t}\left(1+\frac{2t}{n+1}\right)\\
    &\geq {{n-1} \choose t}\frac{(n-2t+3)t}{2(n-t)}
    -{{(n+1)/2} \choose t}\left(2-\frac{2}{n+1}\right)\\
    &> {{n-1} \choose t}\frac{(n-2t+3)t}{2(n-t)}
    -2{{(n+1)/2} \choose t}.
      \end{align*}

The second to last inequality follows from the fact that $$1+\frac{2t}{n+1}\leq 1+\frac{2((n-1)/2)}{n+1} = 2-\frac{2}{n+1},$$
since $t\leq \frac{n-1}{2}$. 

We will show that ${{n-1} \choose t}\frac{(n-2t+3)t}{2(n-t)}
    -2{{(n+1)/2} \choose t}> 0$ by proving the equivalent statement
\begin{equation*}
   \frac{{{n-1} \choose t}}{2{{(n+1)/2} \choose t}}\frac{(n-2t+3)t}{2(n-t)}>1.
\end{equation*}

Expanding binomial coefficients, we obtain
\begin{align}
  &\frac{1}{2}\cdot\frac{(n-1)!}{(n-t-1)!}\cdot\frac{((n+1)/2-t)!}{((n+1)/2)!}\frac{(n-2t+3)t}{2(n-t)}\nonumber\\
  &>\frac{\cancel{(n-t)}(n-t+1)\cdots(n-3)(n-2)(n-1)\cdot \cancel{(n-2t+3)}t}{2\cancel{((n+1)/2-t+1)}((n+1)/2-t+2)\cdots ((n+1)/2-2)((n+1)/2-1)((n+1)/2)\cdot \cancel{2(n-t)}}\nonumber\\
  &= \left[\frac{(n-t+1)}{(n+1)/2-t+2}\cdots\frac{n-3}{(n+1)/2-2}\cdot \frac{n-2}{(n+1)/2-1}\right]\cdot \frac{t(n-1)}{n+1}\label{eq:brackets}
  \end{align}
Notice that $$\frac{n-t+i}{(n+1)/2-t+i+1}=\frac{(n+1)/2-t+i+1+(n-1)/2-1}{(n+1)/2-t+i+1}=1+\frac{(n-1)/2-1}{(n+1)/2-t+i+1}>1$$ for any $i\geq 0$ since $3\leq t\leq (n-1)/2$ and $n\geq 5$. Thus, the term in brackets of \cref{eq:brackets} above is strictly greater than 1, and
\begin{align*}
  &\frac{{{n-1} \choose t}}{2{{(n+1)/2} \choose t}}\frac{(n-2t+3)t}{2(n-t)} >\frac{t(n-1)}{n+1}\ge \frac{2(n-1)}{n+1}=2-\frac{4}{n+1}>1.  \qedhere
  \end{align*}
\end{case}
\end{proofc}
\end{proof}

\section*{Acknowledgments}

The third author thanks Aleyah Dawkins for helpful discussions.

\bibliographystyle{plainurl}
\bibliography{references}

\end{document}

%% file: basics.tex
\newcommand{\Z}{\mathbb{Z}}

\DeclarePairedDelimiter\ceil{\lceil}{\rceil}
\DeclarePairedDelimiter\floor{\lfloor}{\rfloor}

\DeclarePairedDelimiter\set{\{}{\}}

\DeclarePairedDelimiterX\setof[2]{\{}{\}}{#1\,:\,#2}